\documentclass{smfart}

\usepackage[utf8]{inputenc}
\usepackage{amssymb}
\usepackage{amsmath}
\usepackage{amsthm}
\usepackage[shortlabels]{enumitem}
\usepackage{xcolor}
\usepackage{graphicx}

\newtheorem{theorem}{Theorem}
\newtheorem{lemma}{Lemma}
\newtheorem{definition}{Definition}

\newtheorem{proposition}{Proposition}

\newtheorem{corollary}{Corollary}

\newcommand{\secref}[1]{Section~\ref{#1}}
\newcommand{\thmref}[1]{Theorem~\ref{#1}}
\newcommand{\lemref}[1]{Lemma~\ref{#1}}

\newcommand{\propref}[1]{Proposition~\ref{#1}}

\newcommand{\defref}[1]{Definition~(\ref{#1})}

\title{One dimensional non-Hausdorff manifolds and foliations of the plane}
\author{Andr\'e Haefliger and Georges Reeb\\
\footnotesize{(translated by Gangotryi Sorcar)}}

\date{6 December 1956}

\begin{document}

\maketitle

\section*{Foreword from the translator}

The original article titled ``Vari\'et\'es (non séparées) à une dimension et structures feuilletées du plan" was published in 1957 in French in \emph{L'Enseignement Mathématique}. It establishes a beautiful connection between foliations of the plane and non-Hausdorff $1$-dimensional manifolds arising naturally as leaf spaces of the foliations. Since its appearance, this theory has paved the way for several results concerning dynamical systems and foliations of the plane and $2$-manifolds. For example, in $1985$ Neumann and Goel studied unstable dynamical systems of $2$-manifolds. In $2005$, LeRoux studied homeomorphisms of the plane that are fixed point free and orientation preserving (also known as Brouwer homeomorphisms or free mappings). In $2012$, LeRoux, O'Farrell, Roginskaya, and Short developed necessary and sufficient conditions for certain plane homeomorphisms (Reeb homeomorphisms) to be embeddable within a flow. More recently, in $2021$, Groisman and Nitecki have proved a characterization of non-trivial foliations of the plane and used it to study conjugacy classes of Anosov diffeomorphisms of the plane that are time-$1$ maps of flows. Sushil Bhunia and the translator are currently working on classifying conjugacy classes of Brouwer homeomorphisms embeddable in flows.
Haefliger and Reeb's article inspires many results in the theory of foliation of $3$-manifolds as well: we refer the interested reader to Danny Calegari's book on the topic. 

Haefliger and Reeb's article also has been applied to areas outside topological dynamics. For example, in $2008$ Baillif and Gabard showed that connected homogeneous manifolds need not be Hausdorff; in their result, the non-Hausdorff $1$-dimensional manifolds from this article proved to be useful. This article has been referenced in $43$ papers to-date and has helped build many nice results. Only very few of these results are mentioned here due entirely to the translator's inability to convey the remainder in a succinct manner. In the literature, the main theorem of this article has been commonly referred to as Haefliger-Reeb theory or ``a classical result by Haefliger and Reeb." However, to the best of our knowledge, no English version of this article is currently available. We have kept the article entirely unchanged except to add figures and correct a few typographical errors. We hope that this translation will be useful for the broader audience.

\vspace{1cm}





\section*{Introduction}
The concept of a Hausdorff\footnote{Usually we say topological manifold instead of Hausdorff topological manifold.} topological manifold is fundamental to a number of branches of mathematics. Recall its definition:
\vspace{2mm}

A Hausdorff $n$-dimensional topological manifold is a topological space $V^n$ such that: 
\begin{enumerate}
    \item Any point of $V^n$ has an open neighborhood homeomorphic to $\mathbb{R}^n$.
    \item $V^n$ is Hausdorff, i.e. any two points of $V^n$ have open neighborhoods that do not intersect. 
\end{enumerate}

The convenience of condition $2$ appears in the study of certain properties of differential geometry and topology. 

For example, a Hausdorff manifold $V^n$ with a countable basis is metrizable. It follows that it is homeomorphic to a subspace of higher dimensional Euclidean space. Any Hausdorff manifold of one dimension with countable basis is homeomorphic to $\mathbb{R}$ or $\mathbb{S}^1$. 

However, it seems useful to also study topological manifolds that are non-Hausdorff. These spaces arise naturally while answering several questions\footnote[2]{For example, a sheaf defined on a Hausdorff manifold is generally endowed with a non-Hausdorff manifold structure.}. The purpose of this article is to show how the study of one-dimensional manifolds (usually non-Hausdorff) allows us to find several properties of foliations of $\mathbb{R}^2$.

The first section is devoted to the study of non-Hausdorff manifolds (in particular those of one dimension). After providing some definitions and examples \eqref{1.1}, we establish some properties of $1$-dimensional simply connected manifolds \eqref{1.2} and differentiable structures that can be defined on them (1.3). These properties are applied in the second section. 

The foliations of the plane were studied by Poincar\'e and many others. The basic definitions and main results are collected in \eqref{2.1}. The Theorems \ref{kaplan}, \ref{kamke}, and \ref{wazewsky}
are due to Kaplan, Kamke, and Wazewsky, become particularly clear in our opinion if we start from the following fundamental remark: The space of leaves of a foliation of $\mathbb{R}^2$ is a $1$-dimensional manifold (in general non-Hausdorff) \eqref{2.2}. These theorems are presented in \eqref{2.3}.

\section{Properties of $1$-dimensional manifolds}\label{1}

\subsection{Definitions and examples}\label{1.1}

\begin{definition}
An $n$-dimensional topological manifold $V_n$ is a topological space such that any point $p$ admits an open neighborhood homeomorphism to $\mathbb{R}^n$.
\end{definition}

A homeomorphism $h: \mathbb{R}^n\rightarrow U \subset V_n$ where $U$ is the image of $h$ is called a chart of $\mathbb{R}^n$ in $V_n$. The transition map associated with two charts $h_i$ and $h_j$ of $\mathbb{R}^n$ in $V_n$ with respective images $U_i$ and $U_j$ is the homeomorphism $h_j^{-1}h_i$\footnote[3]{If $f$ is a function from a subset $A$ of a set $E$ to a set $E'$ and $f'$ is a function from a subset $A'$ of $E'$ to a set $E''$, we will denote by $f'f$ the function $x\mapsto f'[f(x)]$ of the subset of $E$ formed by points $x$ such that $f(x)\in A'$.}: 

$$h_j^{-1}h_i : h_i^{-1}(U_i\cap U_j)\rightarrow h_j^{-1}(U_i \cap U_j)$$

According to previous definition, there is always a set of charts whose images cover $V_n$. Such a set will be called an atlas of $\mathbb{R}^n$ on $V_n$\footnote[4]{We are using the terminology and definition of C. Ehresmann.}.

Hausdorff topological manifolds often have additional structures: orientation, differentiable structure, complex structure. These notions (as well as, for example, tangent vectors and tensor products) are defined without appealing to the Hausdorff axiom, therefore they translate immediately to non-Hausdorff manifolds.

\vspace{2mm}
\noindent
A general construction process:

\begin{definition}
We say that a continuous function $p:E\rightarrow E'$ ``spreads $E$ over $E'$" if every point $x\in E$ admits an open neighborhood $U$ such that $p|_U$ is a homeomorphism to an open set in $E'$.
\end{definition}

\begin{proposition}\label{prop1}
Let $V_n$ be an $n$-dimensional manifold and let $\rho$ be an equivalence relation in $V_n$ for which every point $x \in V_n$ has a neighborhood $U_x$ such that no two points in $U_x$ are $\rho$-equivalent. Then the topological quotient space $V_n'=V_n/\rho$ is an $n$-dimensional topological manifold, the quotient map $p:V_n\rightarrow V_n/\rho$ is a local homeomorphism, and $V_n$ is an \'etale space over $V_n/\rho$ (or equivalently $p$ spreads $V_n$ over $V_n/\rho$).
\end{proposition}

Indeed, if $U$ is an open neighborhood of $x$ such that a restriction of $\rho$ to $U$ is the identity, the restriction of $p$ to $U$ is a homeomorphism of $U$ to an open set $U'$ in $V_n'=V_n/\rho$; since each point of $V_n$ admits a neighborhood homeomorphic to $\mathbb{R}^n$ so will it be for each point of $V_n'$.

Note that even if $V_n$ is Hausdorff, $V_n'$ may not be, as the examples ahead will show.

The procedure of construction of the manifold given in Proposition $1$ is general in the following sense: 
Given an $n$-dimensional manifold $V_n$, there is a family of charts $h_i:\mathbb{R}^n \rightarrow V_n$ whose images cover $V_n$ (as $i$ takes all values within some index set $I$); the manifold $V_n$ is therefore isomorphic to a quotient of the space $I\times \mathbb{R}^n$ by the equivalence relation associated to the function from $I \times \mathbb{R}^n$ to $V_n$ given by $(i, x) \rightarrow h_i(x)$.

\begin{definition}
A point $x$ in a manifold $V_n$ is called a branch point if there is a point $z \in V_n$ ($z\ne x$) which is not ``separate" from $x$, i.e. any neighborhood of $x$ has non-empty intersection with all neighborhoods of $z$.
\end{definition}

Note: The relation ``$x$ is not separate from $z$" is reflexive and symmetric, but not transitive in general (see Example $3$ below).\\

\noindent{\emph{Examples.}}\\

Proposition \ref{prop1} allows us to construct an array of examples of manifolds. Now, let us study some examples, limiting ourselves to $1$-dimensional manifolds that are second countable.\\

EXAMPLE $1$. Let $R_1$ and $R_2$ be two copies of the real line $\mathbb{R}$ and let $\Sigma$ be the disjoint union of $R_1$ and $R_2$. Consider an open set $\Omega$ in $\mathbb{R}$ and the equivalence relation $\rho$ in $\Sigma$ that identifies points of $R_1$ and $R_2$ that have the same coordinate $t\in \Omega$ and is reduced to identity for the other points, and therefore satisfies the conditions of property of Proposition \ref{prop1}. Passing to the quotient we get a $1$-dimensional manifold (not necessarily Hausdorff). The branch points are those who have coordinates same as boundary points of $\Omega$\footnote[5]{The coordinate of point $x\in \Sigma/\rho$ is the coordinates of points of $\Sigma$ that project onto $x$ by the canonical quotient map $\Sigma \rightarrow \Sigma/\rho$.}.

\begin{center}
    \includegraphics[width=90mm,scale=0.8]{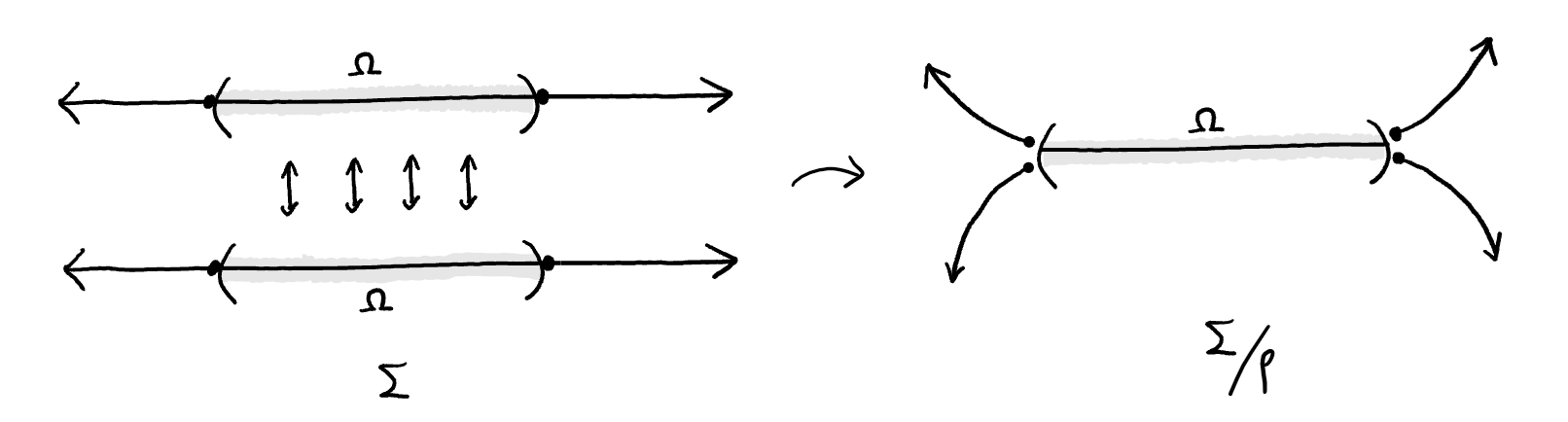}
\end{center}

It is useful to mention some particular cases. 

\begin{enumerate}[a)]
    \item The simple branching: Here $\Omega=(-\infty,0)$. The branch points have coordinate $0$.
    \item The lasso: Here $\Omega=(-\infty, -1)\cup (0,\infty)$. The branch points have coordinates $-1$ and $0$.
    \item The strangled lasso: Here $\Omega=\mathbb{R}\setminus \{0\}$. The branch points have coordinate $0$.
    \item $\Omega$ is the complement of the Cantor's perfect set. Here the collection of branch points is uncountable. 
\end{enumerate}

\begin{center}
    \includegraphics[width=90mm,scale=0.8]{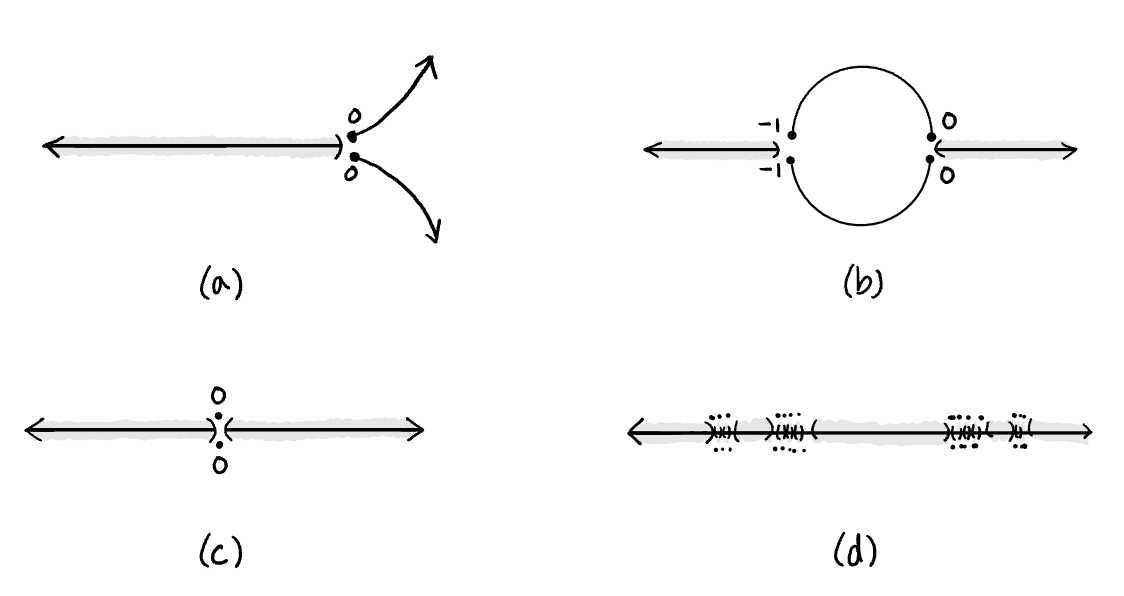}
\end{center}

\vspace{.4cm}

EXAMPLE $2$. The loop: Let $\rho$ be the equivalence relation on $\mathbb{R}$ which identifies the points with coordinates $t$ and $-t$ for $\lvert t \rvert>1$\footnote[6]{The original paper says $\lvert t \rvert<1$ which we have interpreted as a typo since in that case the point $0$ will fail to satisfy the condition of Proposition \ref{prop1}} and is reduced to identity for the other points. The quotient space is the loop and the branch points are $1$ and $-1$.\\

\begin{center}
    \includegraphics[width=50mm,scale=0.55]{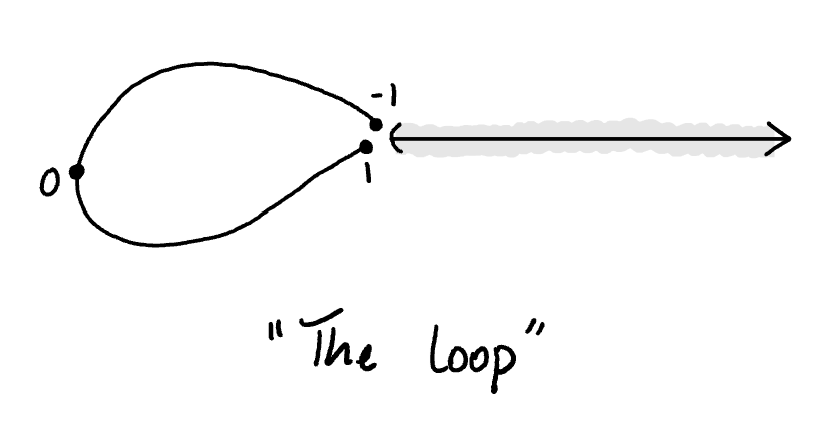}
\end{center}

EXAMPLE $3$. The star: Let $\Sigma$ be the topological sum of $R_1, R_2, \ldots, R_n$, $n$ copies of $\mathbb{R}$. Let $\rho$ be the equivalence relation in $\Sigma$ that identifies each point with coordinate $t>0$ in $R_i$ with the point with coordinate $-t$ in $R_{i+1}$ ($1\le i \le n$; assume $R_{n+1}=R_1$). The quotient space is a $1$-dimensional manifold that looks like a star with $n$ branches. The branch points here are the points with coordinate $0$. Two such points in $R_i$ and $R_j$ are separate if and only if $i-j\ne 1$. We could also consider a star with infinitely many branches.\\

\begin{center}
    \includegraphics[width=100mm,scale=1.9]{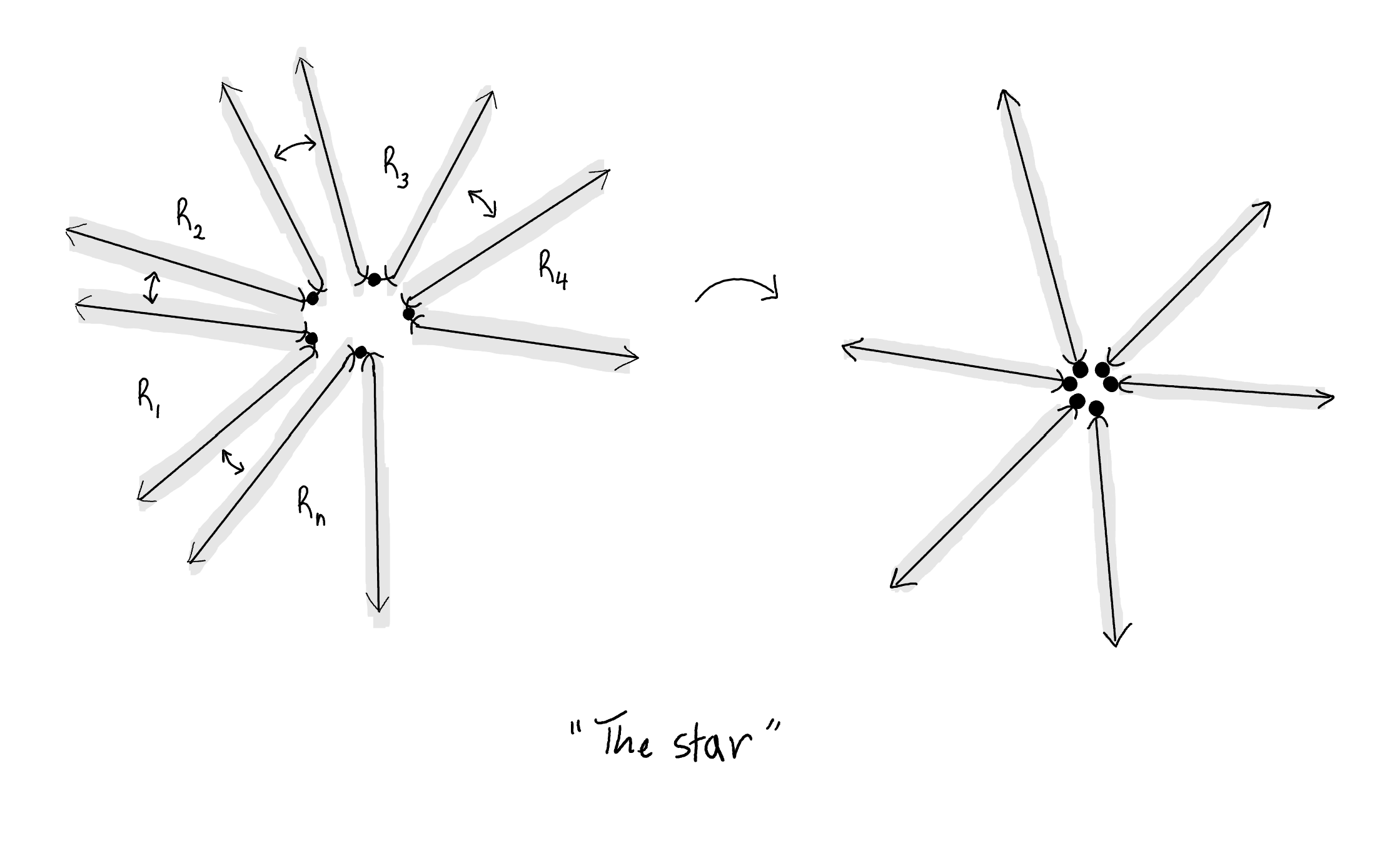}
\end{center}

EXAMPLE $4$: The plume: Example $1 a)$ shows how it is possible to ``graft" a copy of $\mathbb{R}$ as a simple branching at a point $t=0$; we can obviously graft one such branching at a point anywhere in $\mathbb{R}$. If we simultaneously graft a simple branch at all points in $\mathbb{R}$ with rational coordinates -- it is not productive to describe this process in detail -- one obtains a manifold of $1$-dimension that we call ``the plume". The straight line $\mathbb{R}$ is the stem on which the branches are grafted. Here the branch points form a dense set in $\mathbb{R}$.

\begin{center}
    \includegraphics[width=100mm,scale=1.9]{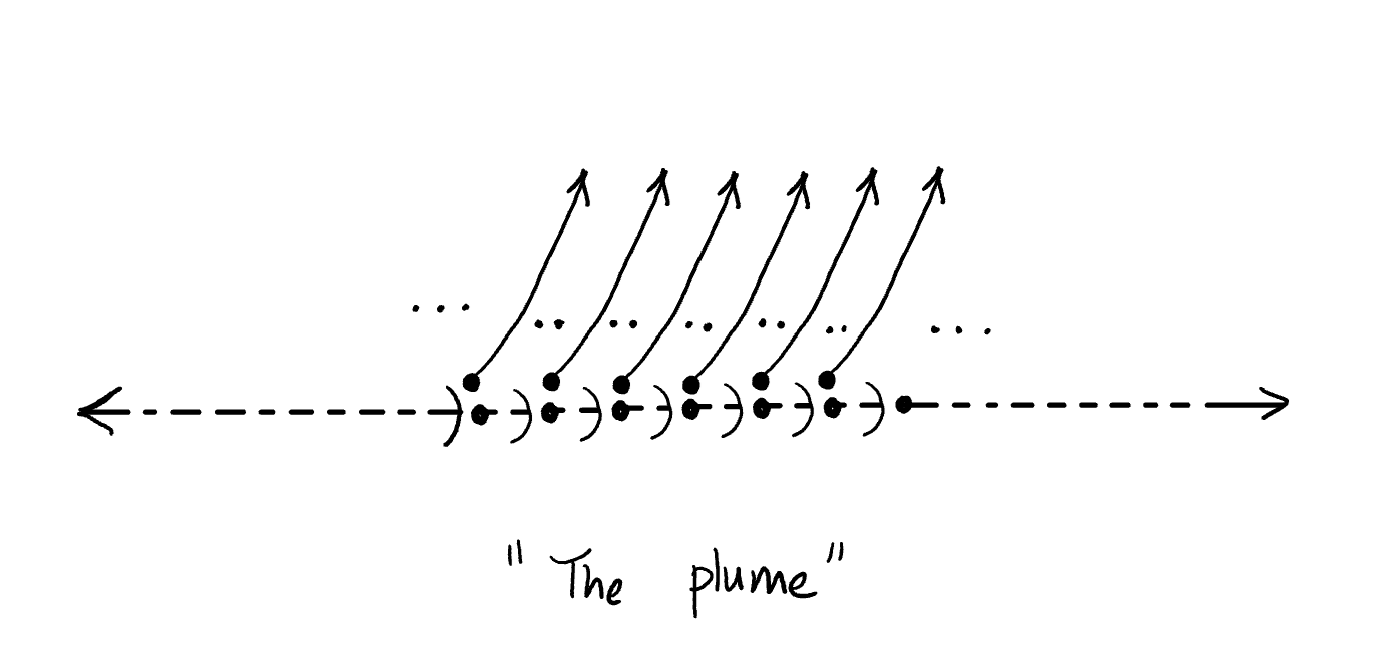}
\end{center}

By grafting a branch at every point of $\mathbb{R}$ we would obtain a manifold that does not have a countable basis.\\

EXAMPLE $5$: The compound plume: If in a plume we replaced each branch with a new plume, we get a new manifold of $1$-dimension that we call the double plume. In a double plume you can replace each branch with a simple plume, thus obtaining a triple plume. Repeating this process $n$ times ($n$ being a natural number) we get the $n$-plume. We can carry out a countable sequence of these operations to get the $1$-dimensional manifold that is called the compound plume. The compound plume has the following remarkable property: The branch point set is countable and everywhere dense in the obtained space. 

\begin{center}
    \includegraphics[width=100mm,scale=1.9]{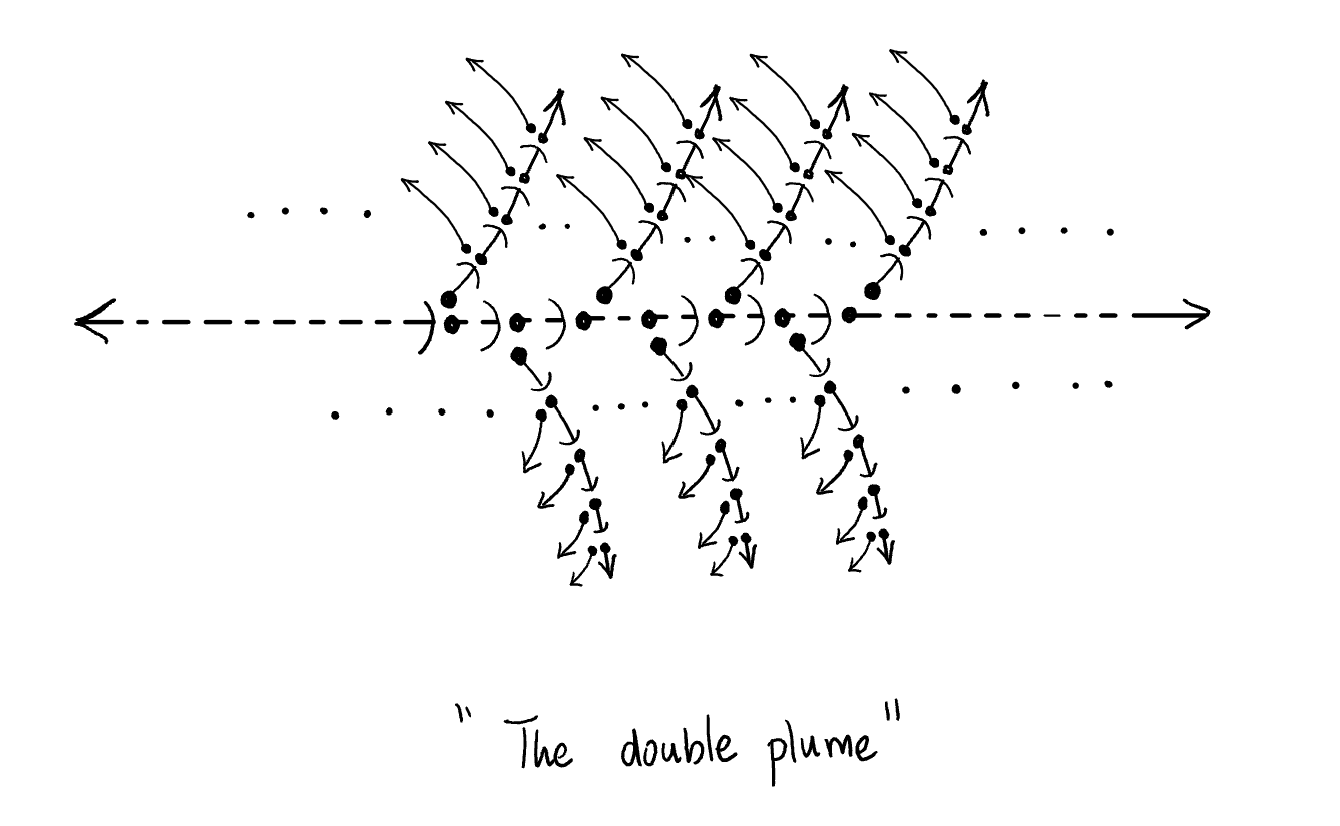}
\end{center}

These examples show the great diversity of connected, second countable, $1$-dimensional manifolds. A topological classification of these spaces already seems quite complicated. 

Let us give another example which will give a faint idea of the complexity of non-Hausdorff manifolds of dimension greater than $1$. Let $E_1$ and $E_2$ be two copies of $\mathbb{R}^2$ with the polar coordinate system $(r, \omega)$. Let $\rho$ be the equivalence relation on the space $\Sigma=E_1\sqcup E_2$ that identifies any point of $E_1$ with coordinate $(r, \omega)$, $r<1$ with the point of $E_2$ with coordinate $(r, \frac{\omega}{1-r})$ and is identity for other points. The quotient space $\Sigma/\rho$ is a $2$-dimensional manifold; the image in $\Sigma/\rho$ of any point of $E_1$ with coordinate $(1, \omega_0)$ is not Hausdorff from any point of the image in $\Sigma/\rho$ of the circle $r=1$ of $E_2$. 

\subsection{Simply connected $1$-dimensional manifolds}\label{1.2}

Recall the following definitions:

\begin{definition}
The pair $(\widetilde{V}, p)$ of a topological space $\widetilde{V}$ and a continuous function $p:\widetilde{V}\rightarrow V$, where $V$ is another topological space, is called a covering of $V$ if every  point of $V$ has an open neighborhood $U$ such that $p^{-1}(U)$ admits a partition in open subsets $U_i$ such that the restriction of $p$ on each $U_i$ is a homeomorphism onto $U$.
\end{definition}

A topological space $V$ will be called simply connected if it is connected and if, for any connected covering $(\widetilde{V}, p)$ of $V$, the projection $p$ is a homeomorphism of $\widetilde{V}$ to $V$.

\begin{definition}
An $n$-dimensional manifold $V_n$ is said to be orientable if there exists an atlas $A$ of $\mathbb{R}^n$ on $V_n$ such that any change of charts associated to two charts in $A$ is an orientation preserving homeomorphism from an open set of $\mathbb{R}^n$ to an open set of $\mathbb{R}^n$.
\end{definition}

The manifolds constructed in examples $2$ and $3$ (for $n$ odd) are not orientable.

If a $1$-dimensional manifold $V$ can be ``spread" over (for which there exists a local homeomorphism $V\to \mathbb{R}$) the real line $\mathbb{R}$, then $V$ is necessarily orientable. On the other hand, an orientable $1$-dimensional manifold does not always spread over $\mathbb{R}$: just consider the case of the circle. However, it is possible to show the following proposition which will be useful later. 

\begin{proposition}\label{localhomeo}
Suppose $V$ is a $1$-dimensional simply connected manifold with countable basis. Then there exists a continuous function $f:V\rightarrow \mathbb{R}$ which is a local homeomorphism.
\end{proposition}

In the proof, we will use the following lemma:

\begin{lemma}\label{lem1}
If $V$ is a $1$-dimensional simply connected manifold, then for any $x\in V$ $V\setminus\{x\}$ has two connected components.
\end{lemma}

Indeed, let $U$ be a neighborhood of $x$ homeomorphic to an interval; $U\setminus\{x\}$ has two connected components $U_+$ and $U_-$. Consider now two copies $V'$ and $V''$ of $V\setminus\{x\}$ and let $U_+'$, $U_-'$ and  $U_+''$, $U_-''$ be the copies of $U_+$ and $U_-$ in $V'$ and $V''$ respectively. Complete the space $V' \sqcup V''$ by two points $x'$ and $x''$ admitting respectively the neighborhoods $U'$ and $U''$ such that $U'\cap V'=U_+'$, $U'\cap V''=U_-''$ and $U''\cap V'=U_-'$, $U''\cap V''=U_+''$ to obtain a space $\overline{V}$ that, with the canonical projection $p$ on $V$ (in particular $p(x')=p(x'')=x)$, is a $2$-sheeted cover of $V$\footnote[7]{The proof can be completed by noting that if $V\setminus x$ is connected, then $\overline V$ is connected, hence $V$ is not simply connected.}.

\begin{center}
    \includegraphics[width=130mm,scale=3.5]{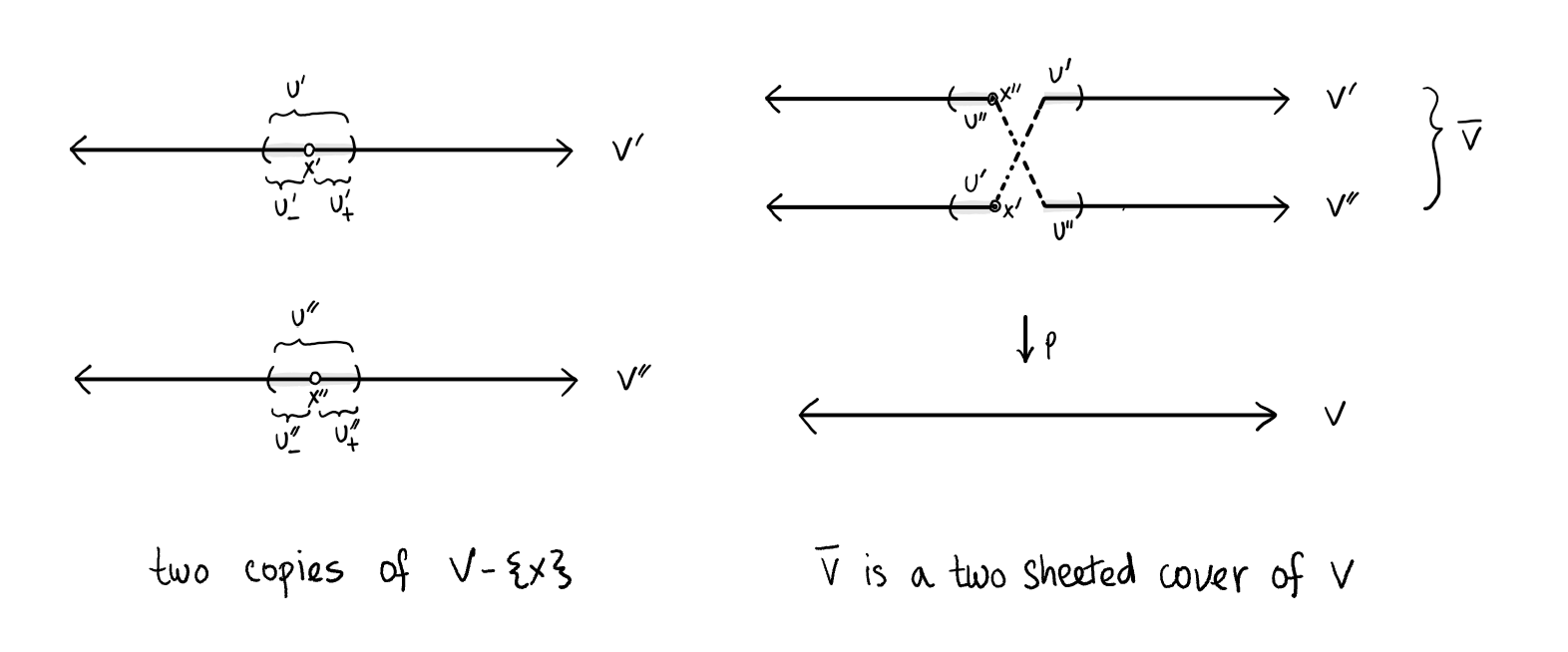}
\end{center}

Conversely, one can show that if $V$ is a $1$-dimensional connected manifold with the property that the respective complement of every point in $V$ is disconnected, then $V$ is simply connected.

Let us now proceed to the proof of the proposition. It is possible to find a countable family of charts $h_i$   $(i=1, 2,\ldots)$ from $\mathbb{R}$ to $V$ whose images $O_i$ cover $V$. Since $V$ is connected we can assume that the numbering of $O_i$ is made such that the union $\Omega_n=\bigcup_{i=1,2,\ldots,n}O_i$ is connected for all integers $n$. Now we will reason by induction. Suppose we define on $\Omega_n$ a continuous function $f_n$ that spreads $\Omega_n$ over the interval $(-n,n)$. We will show that $f_n$ can be further extended to the function $f_{n+1}$ that spreads $\Omega_{n+1}$ over the interval $(-n-1, n+1)$.

It follows from the lemma that $\Omega_n \cap O_{n+1}$ is connected because if this were not the case then we could find a point $x$ such that $V\setminus\{x\}$ is connected. So $h_{n+1}^{-1}(\Omega_n \cap O_{n+1})=I$ is an open interval of $\mathbb{R}$ on which the function $f_nh_{n+1}$ is continuous, strictly monotonic, and such that $\lvert f_nh_{n+1}\rvert<n$. The function $f_nh_{n+1}$ can be extended to a continuous function $\varphi$ on $\mathbb{R}$, strictly monotonic and for which $\lvert \varphi \rvert<n+1$. The function $\varphi h_{n+1}^{-1}$\footnote[8]{The original paper writes $\varphi h^{-1}$ which is a typo.} defined on $O_{n+1}$ and the function $f_n$ are identical on $\Omega_n\cap O_{n+1}$, so amalgamating them defines the desired extension $f_{n+1}$. 

\subsection{Manifolds with differentiable structures}\label{1.3}

\begin{definition}
A $C^r$-differentiable structure for $r$ a positive integer or $\infty$ (respectively an analytic structure) on an $n$-dimensional manifold $V_n$ is defined by the data of an atlas $A$ of $\mathbb{R}^n$ on $V_n$ such that, for any two charts $h_i, h_j \in A$, the transition map $h_j^{-1}h_i$ is a homeomorphism that is $r$ times continuously differentiable (resp. analytic) (i.e., $h_j^{-1}h_i$ is a $C^r$ (resp. analytic) homeomorphism) from an open set of $\mathbb{R}^n$ to an open set of $\mathbb{R}^n$.
\end{definition}

An $r$-differentiable function on $V_n$ is a function $f:V_n\rightarrow \mathbb{R}$ such that for any chart $h_i\in A$, the function $fh_i$ is a $C^r$ function on $\mathbb{R}^n$. An $r$-differentiable function or $C^r$ function on $V_n$ is called ``rank $1$'' at a point $x\in V_n$ if, for one chart $h_i$ whose image contains $x$, the function $fh_i$ has at least one nonzero partial derivative at the point $h_i^{-1}(x)$. This definition is obviously independent of the choice of chart $h_i\in A$. 

We similarly define the notion of $C^r$ functions from a $C^r$-differentiable manifold $V_n$ to a $C^r$-differentiable manifold $V_m$. 

A chart $f:\mathbb{R}^n\rightarrow V_n$ is called ``compatible'' with atlas $A$ if, for all $h\in A$, the transition maps $f^{-1}h$ and $h^{-1}f$ are $C^r$ (or analytic) homeomorphisms from open sets of $\mathbb{R}^n$ to $\mathbb{R}^n$. The set of all compatible charts with $A$ form a maximal atlas $\overline{A}$ generated by $A$. Two sub-atlases of $\overline{A}$ define the same $C^r$-differentiable (or analytic) structure on $V_n$.

Property $1$: The simple branching (Example $1 a)$) is a $1$-dimensional manifold admitting two non-isomorphic (i.e. distinct) $C^\infty$-manifold structures. (A similar argument applies to the strangled lasso.)

In other words, the simple branching can be given two $C^\infty$ differentiable structures such that there is no $C^\infty$ homeomorphism (with a $C^{\infty}$ inverse) from $V$ with the first structure to $V$ with the second structure. 

Using the notations of the Example $1$, $R_1$ and $R_2$ (identified to their images in $V$) are two open sets forming a cover of $V$; let the point of $R_1$ (respectively $R_2$) of coordinate $t$ be designated $t_1$ (respectively $t_2$). One can define a $C^\infty$-differentiable structure on $V$ by giving two maps $h_1$ and $h_2$ on $\mathbb{R}$ to $R_1$ and $R_2$ respectively such that $h_1^{-1}h_2$ and $h_2^{-1}h_1$ are $C^\infty$-homeomorphisms of the half line $(-\infty, 0)$. 

First structure: $h_1(t)=t_1$ and $h_2(t)=t_2$.

Second structure: $h_1(t)=t_1$ and $h_2(t)=t_2^3$. 

For the first structure, the function which takes the value $t$ at the point $t_1$ or $t_2$ is a $C^\infty$-function of rank $1$ everywhere on $V$. On the other hand, for the second structure, any $C^\infty$-function $f$ on $V$ is rank $0$ at the point $t=0$. Indeed, let $f_1=fh_1$ and $f_2=fh_2$: from $f_1(t)=f_2h_2^{-1}h_1(t)$ we get $\left ( \frac{d}{dt}f_1(t)\right )_{t=0}=0$. Alternatively, one can use $f_2=f_1h_1^{-1}h_2(t)$ to obtain $\left ( \frac{d}{dt}f_2(t)\right )_{t=0}=0$ and we thank the referee for pointing this out. This establishes Property $1$\footnote[9]{In a recent article Milnor built two distinct differentiable structures on $\mathbb{S}^7$. Milnor's result is global, and on the contrary, here we have a local property relative to an arbitrary open neighborhood of a pair of points that cannot be separated.}.

Note that the two structures defined above are even analytic. It is clear that we can produce multiple examples. We understand how to easily use the previous example to construct a $C^\infty$-differentiable structure on a compound plume which highlights the following property:

Property $2$: There exists $1$-dimensional manifolds (for example, the compound plume) with $C^\infty$-differentiable structures such that continuously differentiable functions on these manifolds are constant functions. 

The pathological properties highlighted above lead to a more strict notion of ``differentiable structure'' for which the previous properties are no longer valid. 

\begin{definition}
A $C^r$-differentiable structure on a manifold $V_n$ is called \emph{regular} if for all $C^r$-differentiable functions $f$ defined on a neighborhood of $x\in V_n$, there is a $C^r$-differentiable function $f'$ defined on $V_n$ such that $f$ and $f'$ coincide on a neighborhood of $x$. 
\end{definition}

All $C^r$-differentiable structures on Hausdorff manifolds are regular. The second structure on the simple branching is not regular. 

In the proof of Proposition $2$ we will use the lemma below:

\begin{lemma}
Let $V$ be a manifold with a regular $C^r$-differentiable structure. If a real-valued $C^r$ function on $V$ is of rank $1$ at a point $x\in V$, it is also of rank $1$ at any point $y$ not separated from $x$.
\end{lemma}

To simplify the notation we will prove this lemma in the case where $V$ is a $1$-dimensional manifold. Let $h_1$ and $h_2$ be two charts of $\mathbb{R}$ in $V$ such that $h_1(0)=x$ and $h_2(0)=y$. The function $h=h_2^{-1}h_1$ (as well as its inverse) is a $C^r$-homeomorphism from an open set $U_1$ of $\mathbb{R}$ to an open set $U_2$ of $\mathbb{R}$, the origin belonging to the closure of $U_1$ and $U_2$; the functions $f_1=fh_1$ and $f_2=fh_2$ are $C^r$-differentiable on $\mathbb{R}$ and $f_1'(0)\ne 0$. Suppose $t_1, t_2, \ldots, t_n, \ldots$ is a sequence of points in $U_1$ such that $\displaystyle\lim_{n\to\infty}t_n=0$; by setting $\overline{t}_n=h(t_n)$ we also have $\displaystyle\lim_{n\to\infty}\overline{t}_n=0$. As $f_1'(t_n)=f_2'(\overline{t}_n)h'(t_n)$ for all $n$, if $\displaystyle\lim_{n\to\infty}f_2'(\overline{t}_n)=f_2'(0)$ is $0$, then $\displaystyle\lim_{n\to\infty}h'(t_n)=0$ would be infinite; but then if $g$ is a $C^r$ function on $V$ of rank $1$ at $y$ (such a function always exists by virtue of the assumption of regularity), with the corresponding notations, $\displaystyle\lim_{n\to\infty}g_1'(t_n)$ would be infinite, which is impossible. So $f_2'(0)\ne0$. 

Proposition $2$ can be restated in the following way:
\vspace{2mm}

\noindent
\textbf{Proposition $2$:} \emph{On any $1$-dimensional regular $C^r$-differentiable manifold $V$ that is simply connected and has a countable basis there is a $C^r$-differentiable function with rank $1$ everywhere.} 

In other words, the manifold $V$ can be spread over $\mathbb{R}$ by a $C^r$ function that is rank $1$ everywhere. 

We assume the following known lemma:

\begin{lemma}
Let $f(t)$ be a $C^r$ function defined on an open set $I'$ in $\mathbb{R}$ and whose derivative is nonzero at any point in a closed interval $I\subset I'$; the restriction of $f$ to $I$ can be extended to a $C^r$ function on $\mathbb{R}$ whose derivative is nonzero everywhere on $\mathbb{R}$.
\end{lemma}

Suppose $A$ is the atlas of $\mathbb{R}$ on $V$ that defines the $C^r$-differentiable structure on $V$. The general setup of the proof is like that of Proposition $1$ of $1.2$. Let us take the same notations this time assuming that each $O_i$ is the image of the interval $(-1,1)$ by a homeomorphism $h_i$ which extends along a homeomorphism $\widetilde{h}_i\in A$ of $\mathbb{R}$ in $V$. We suppose a $C^r$ function $f_n$ of everywhere rank $1$ is defined on $\Omega_n$ and that for all $O_i$, $1\le i \le n$, the function $f_nh_i$ which is defined on $I=(-1,1)$ may be extended to a $C^r$ function on $\mathbb{R}$ with rank $1$ everywhere. We will show that $f_n$ can be extended to a $C^r$ function $f_{n+1}$ on $\Omega_{n+1}$ with the same properties.

As $\Omega_{n}\cap O_{n+1}$ is connected, $f_nh_{n+1}$ is a $C^r$ function over an interval $(t_0,t_1)$ contained in $I$ with  nonzero derivative. The sets $h_{n+1}^{-1}(O_i)$, $1\le i \le n$ are open intervals that cover $(t_0,t_1)$. Suppose $O_k$ is an open set such that $h_{n+1}^{-1}(O_k)$ is an interval of the form $(t_0,t_2)$ where $t_2\le t_1$. By the induction hypothesis, the function $f_nh_k$ can be extended on $\mathbb{R}$ to a function $\widehat{f}_k$ with rank $1$ everywhere. Let $t_0'$ be the point of the interval $[-1,1]$ defined by $t_0'=\displaystyle\lim_{t\to t_0} h_k^{-1}h_{n+1}(t)$ and let $x'=\widetilde{h}_k(t_0')$; the point $x=\tilde{h}_{n+1}(t_0)$ is not Hausdorff from $x'$. There exists a $C^r$ function $g$ on $V$ which coincides with $\widehat{f}_k\widetilde{h}_k^{-1}$ on a neighborhood of $x'$, and as $g$ is of rank $1$ at $x'$, it is also rank $1$ at $x$ (Lemma $2$). The function which is equal to $gh_{n+1}$ on the interval $(-\infty, t_0]$ and to $f_nh_{n+1}$ on $(t_0,t_1)$ is $C^r$ on $(-\infty, t_1)$ because the two functions $g\widetilde{h}_{n+1}$ and $f_nh_{n+1}$ coincide in $(t_0, t_2')$. By repeating the same construction for $t_1$, we get a function defined on $\mathbb{R}$ of rank $1$ on $(t_0, t_1)$ and whose restriction to $(t_0,t_1)$ coincides with $f_nh_{n+1}$; according to Lemma $3$, there exists a $C^r$ function $\widehat{f}_{n+1}$ that extends $f_nh_{n+1}$ that is everywhere rank $1$. The function $f_{n+1}$ we seek is equal to $f_n$ on $\Omega_n$ and to $\widehat{f}_{n+1}h_{n+1}^{-1}$ on $O_n$. 

\begin{corollary}
All differentiable structures on $\mathbb{R}$ are equivalent. 
\end{corollary}

Let $R$ be $\mathbb{R}$ with the usual differentiable structure on $\mathbb{R}$ and $R'$ be $\mathbb{R}$ with a $C^r$-differentiable structure. According to the proposition there exists a $C^r$ function of rank $1$ everywhere from $R'$ to $R$ (up to a suitable homothety). As this function is one-to-one, it is a $C^r$ isomorphism of $R'$ and $R$ (provided with the usual $C^r$-differentiable structure).

\section{Foliations of $\mathbb{R}^2$}

\subsection{Review of definitions and classic properties}\label{2.1}

\begin{definition}
A foliation $(F)$ on a $2$-dimensional manifold $V_2$ is defined by an atlas $A$ from $\mathbb{R}^2$ to $V_2$ such that if $h_i$ and $h_j$ are any two charts in $A$, the transition map $h_{ji}=h_j^{-1}h_i$ is a homeomorphism from open set $U_{ji}$ of $\mathbb{R}^2$ to an open set of $\mathbb{R}^2$ such that, in a neighborhood of any point of $U_{ji}$, it is expressed by equations of the form: \begin{equation}\label{eq1}
    x'=g_{ji}(x,y) \hspace{1cm} y'=k_{ji}(y)
\end{equation}
If the transition maps $h_{ji}$ are $C^r$ (resp. analytic), we will say that the foliation $(F)$ is $C^r$ (resp. analytic)\footnote[10]{For a general definition of foliated manifolds, see \cite{R}}.
\end{definition}

Let's say that a chart $f$ of $\mathbb{R}^2$ in $V_2$ is compatible with $A$ if for any chart $h\in A$, the transition map $f^{-1}h$ is also of the form $(1)$ (and furthermore is $C^r$, resp. analytic, if $(F)$ is $C^r$, resp. analytic). All of the atlases compatible with $A$ form the maximal atlas generated by $A$ and defines on $V$ the foliation $(F)$. In what follows, we will assume $A$ is a maximal atlas. 

Let $U$ be an open set of $V_2$. The set of charts of $A$ whose images lie in $U$ form an atlas of $\mathbb{R}^2$ on $U$ that defines the foliation $(F_U)$ induced from $(F)$ on $U$. 

If $(F)$ and $(F')$ are foliations of $V_2$ and $V_2'$ respectively defined by the maximal atlases $A$ and $A'$, an isomorphism of $(F)$ to $(F')$ is a homeomorphism $\psi: V_2\rightarrow V_2'$ such that $A'=\psi A$ (that is any chart of $A'$ is of the form $\psi h$, where $h\in A$). 

In the image $O_i$ of each chart $h_i\in A$, we define an equivalence relation $\rho_i$ whose classes are the images of $h_i$ restricted on straight lines $y=\text{constant}$. From \eqref{eq1} it follows that for any point $x\in O_i \cap O_j$, the equivalence relations induced by $\rho_i$ and $\rho_j$ coincide in a sufficiently small neighborhood of $x$. Let $\rho$ be the equivalence relation generated by the relations $\rho_i$.

\begin{definition}\label{leaves}
The classes of $\rho$ in $V_2$ are called leaves of the foliation $(F)$. 
\end{definition}

The leaf space, i.e. the quotient space of $V_2$ by the equivalence relation $\rho$ (provided with the quotient topology of that of $V_2$) will play an essential role in the following. 

Note that the equivalence relation $\rho$ is open since it is generated by the $\rho_i$, which are open. Recall that for any vector field $E$ defined on a Hausdorff manifold $V_2$ satisfying the following two conditions:

\begin{enumerate}[(i)]
    \item $E$ is $C^r$-differentiable (or analytic),
    \item $E(z)\ne 0$ for all $z\in V_2$,
\end{enumerate}

there is an associated foliation which is $C^r$-differentiable (or analytic). The leaves are then vector field trajectories of $E$. Reciprocally, to all $C^r$-foliations $(F)$ of a Hausdorff manifold $V_2$ \footnote[11]{The original paper writes $\mathbb{R}^2$ instead of $V_2$. We thank the referee for pointing this out.} we can match a vector field $E$ on $V_2$ whose trajectories are the leaves of $(F)$.

{\em Example:} The curves which are solutions of the differential equation (in polar coordinates $r$ and $\omega$) $\frac{dr}{d\omega}=r(1-r^2)$ are the leaves of an analytic foliation of $\mathbb{R}^2\setminus \{0\}$. The circle $r=1$ is one leaf around which the other leaves curl asymptotically. 

The quotient space of $\mathbb{R}^2\setminus \{0\}$ by the equivalence relation $\rho$ associated with the above foliation admits a partition into an open subspace homeomorphic to two circles and one point whose only neighborhood is the entire space.

\begin{definition}
The pair $(O_i, h_i)$ formed by a chart $h_i\in A$ and its image $O_i$ is called a distinguished open set of the foliation $(F)$.
\end{definition}

Let us state the main results related to the foliations of $\mathbb{R}^2$. \thmref{poincarebendixon} below is classical; its proof is based on Jordan's theorem (a particularly easy application); it basically uses the fact that $\mathbb{R}^2$ is simply connected (or more precisely, that its $\beta_1$ (first Betti number) is even). 

\begin{theorem}[Poincar\'e, Bendixon]\label{poincarebendixon}
  Let $(O_i, h_i)$ be a distinguished open set of a foliation of $\mathbb{R}^2$. The image of $h_i^{-1}$ of the intersection of $O_i$ with any leaf is reduced to the empty set or to a line $y=\mathrm{constant}$.
\end{theorem}

\begin{theorem}[Kaplan]\label{kaplan}
  To any foliation of $\mathbb{R}^2$ we can associate a real valued function $\psi:\mathbb{R}^2\rightarrow \mathbb{R}$ with the following properties: 
  \begin{enumerate}[(i)]
      \item $\psi$ is continuous and does not admit a maximum or a minimum. 
      \item $\psi$ is constant on the leaves of the foliation. 
  \end{enumerate}
\end{theorem}

\begin{theorem}[Kamke]\label{kamke}
  Let $(F)$ be a $C^r$-foliation of $\mathbb{R}^2$. Let $\Omega$ be a bounded open set of $\mathbb{R}^2$. There exists a function $\psi:\Omega\rightarrow \mathbb{R}$ satisfying the following properties: 
  \begin{enumerate}[(i)]
      \item $\psi$ is $C^r$ and the gradient of $\psi$ is nonzero at all points of $\Omega$.
      \item $\psi$ is constant on the leaves of the foliation induced by $(F)$ on $\Omega$. 
  \end{enumerate}
\end{theorem}

\begin{theorem}[Wazewsky]\label{wazewsky}
  We can equip $\mathbb{R}^2$ with a $C^\infty$-foliation such that any $C^r$ function on $\mathbb{R}^2$ that is constant on the leaves of the foliation is reduced to a constant function. 
\end{theorem}

We will not reproduce the proof of Theorem \ref{poincarebendixon} but we will show in \ref{2.2} and \ref{2.3} that Theorems \ref{kaplan}, \ref{kamke}, and \ref{wazewsky} are consequences of Theorem \ref{poincarebendixon} and the properties of $1$-dimensional manifolds (as) described in the first part.

Examples of foliations of $\mathbb{R}^2$:
\begin{enumerate}
    \item The lines $y=\text{constant}$ are obviously the leaves of a foliation of $\mathbb{R}^2$.
    \item Let $C$ be a Jordan curve in $\mathbb{R}^2$. The previous foliation induced by $C$ is homeomorphic to an analytic foliation on $\mathbb{R}^2$.
    \item The complement $U$ in $\mathbb{R}^2$ of the set of points with coordinates $x=0$ and $y\ge 0$ is homeomorphic to $\mathbb{R}^2$. The connected components of the level curves of the function $\psi=xy$ are the leaves of an analytic foliation on $U$. 
\end{enumerate}

\subsection{The leaf space of a foliation of $\mathbb{R}^2$}\label{2.2}
The following proposition is an essential consequence of \thmref{poincarebendixon} of \secref{2.1}.

\begin{proposition}\label{leafspacestructure}
Let $(F)$ be a foliation of $\mathbb{R}^2$. The quotient space $V$ of $\mathbb{R}^2$ by the equivalence relation $\rho$ associated with the foliation (see \defref{leaves} of \secref{2.1}) is a $1$-dimensional manifold that is simply connected and has a countable basis. If $(F)$ is a $C^r$-foliation (or analytic foliation) the leaf space $V$ is provided with a canonical $C^r$ (or analytic) manifold structure. 
\end{proposition}

As $\mathbb{R}^2$ is connected and has countable basis, $V$ is also connected and has countable basis. To show that $V$ is a $1$-dimensional manifold it will suffice to show that all points $z$ in $V$ admit an open neighborhood homeomorphic to $\mathbb{R}$. Let $\pi:\mathbb{R}^2\rightarrow V$ be the canonical projection; the leaf $\pi^{-1}(z)$ meets at least one distinguished open set $O_i$. The equivalence relation induced by $\rho$ on $O_i$\footnote[12]{The original paper writes $O^i$ but we think it is a typo and it should be $O_i$.} is the relation $\rho_i$, according to Theorem $1$. So $\pi(O_i)$, which is an open neighborhood of $z$ because $\rho$ is an open equivalence relation, is homeomorphic to $O_i/\rho_i$ and therefore to $\mathbb{R}$.

By virtue of Jordan's theorem, the complement of any leaf (which is a closed subset of $\mathbb{R}^2$) has two connected components; the complement of any point of $V$ therefore also has two connected components. The property is equivalent to the fact that $V$ is simply connected (see \lemref{lem1}).

Let $A$ be an atlas that defines a differentiable foliation on $\mathbb{R}^2$, and let $\overline{h}_i$ be the restriction of the chart $h_i\in A$ on the line $x=0$. The set of charts $\pi\overline{h}_i$ of $\mathbb{R}$ in $V$ is an atlas that defines a differentiable structure on $V$. 

The reader can construct as an exercise the space of leaves of the foliations of $\mathbb{R}^2$ defined in the above examples. 

\subsection{Applications of results of Section \ref{1} to foliations of $\mathbb{R}^2$}\label{2.3}

\subsection*{Proof of \thmref{kaplan}:}

\thmref{kaplan} of \ref{2.1} (Kaplan) is an immediate consequence of \propref{localhomeo} and \propref{leafspacestructure}. Let $f$ be a function that spreads $\mathbb{R}^2/\rho$ over $\mathbb{R}$ (i.e. local homeomorphism); the function $\varphi=f\pi$ is a real valued function on $\mathbb{R}^2$ that satisfies the conditions of \thmref{kaplan}. 

\subsection*{Proof of \thmref{kamke}:} If the theorem is true for bounded open sets that are simply connected in $\mathbb{R}^2$, it will also be true for all bounded open sets $\Omega$: it suffices to consider the open set $C$ bounded by a circle which contains $\Omega$; if $\varphi$ is a function in $C$ which satisfies the condition of the theorem, it will obviously do so for its restriction to $\Omega$. So Kamke's theorem is a consequence of the restated version of \propref{localhomeo} in Section $1.3$ and of the following proposition:

\begin{proposition}
The space of leaves $V'=\Omega/\rho'$ of the foliation $(F')$ on $\Omega$ induced by the foliation $(F)$ is provided with a regular $C^r$-differentiable manifold structure. 
\end{proposition}

Let $\pi$ and $\pi'$ respectively be canonical projections of $\mathbb{R}^2$ on $\mathbb{R}^2/\rho=V$ and $\Omega/\rho'=V'$. The canonical injection $\Omega\rightarrow \mathbb{R}^2$ defines a continuous function $\psi: V'\rightarrow V$ by passing to quotients. Since any distinguished open set of $\Omega$ is a distinguished open set in $\mathbb{R}^2$, the function $\psi$ spreads $V'$ over $V$ (a local homeomorphism) and it is $C^r$ and everywhere rank $1$. Let $x'$ be a point of $V'$ and let $x=\psi(x')$. As $\Omega$ is relatively compact, the intersection of the leaf $F=\pi^{-1}(x)$ with the closure $\overline{\Omega}$ of $\Omega$ is compact. We can then find a distinguished open set $W$ in $\mathbb{R}^2$ sufficiently narrow and elongated to contain $F\cap\overline{\Omega}$ and such that $W\cap \Omega$ is saturated by leaves of $\Omega$. Let $U:=\pi(W)$, an open set containing $x$ homeomorphic to an interval. Let $f'$ be a $C^r$ function defined on an open set $U'$ containing $x' \in V'$ homeomorphic to an interval and sufficiently small that $\psi(U')\subset U$. The function $f=f' \circ \psi|_{U'}^{-1}$ is a $C^r$ function defined on an open set containing $x$. Since $U$ is homeomorphic to an interval, it is possible to construct a $C^r$ function $g$ on $U$ that coincides with $f$ on an open neighborhood of $x$ and that is zero outside a compact set $K$ containing $U$. As $\pi^{-1}(K)\cap W\cap\Omega$ is closed in $\Omega$ and is saturated with leaves of $\Omega$, the function which is equal to $g\pi$ on $W\cap \Omega$ and zero at other points of $\Omega$ is $C^r$, hence passing to quotients we obtain a $C^r$ function $g'$ on $V'$ which coincides with $f'$ in a neighborhood of $x'$.

\subsection*{Proof of \thmref{wazewsky}:} It is generally necessary to assume that $\Omega$ is bounded for Kamke's theorem to be true, as the following example shows. 
Let $O_1$ and $O_2$ be two copies of $\mathbb{R}^2$ in the disjoint union $O_1\sqcup O_2$. Consider the equivalence relation $O_1\ni (x_1,y_1)\sim(x_2,y_2)\in O_2$ for $x_2=x_1+1/y_1$,  $y_2=y_1^3$ for all $y_1<0$, and identity for all points such that $y_1\ge 0$ or $y_2\ge 0$. The quotient space $E$ is homeomorphic to $\mathbb{R}^2$; the two charts $O_1$ and $O_2$ form an atlas that defines an analytic foliation on $E$ where the leaf space is the simple branching with the second differentiable structure defined in Property $1$ of Section $1.3$.

This example helps to understand the construction of the example by Wazewsky. In particular, it's easy to imagine a $C^\infty$-foliation with leaf space a composite plume (Example $5$ of Section $1.1$) equipped with a $C^\infty$-differentiable structure such that any differentiable function on $V$ is reduced to a constant.

\vspace{.5cm}
\begin{center}
   CLASSIFICATION OF THE FOLIATIONS OF $\mathbb{R}^2$.
\end{center}
\vspace{.2cm}

The problem of classification of foliations of $\mathbb{R}^2$ was completely resolved by Kaplan. We will indicate briefly and without proof how our methods also solve this problem. 

It is insufficient to consider only the quotient space $V$ associated with the foliation to characterize the foliation. For that we need to introduce an order relation among the branch points. Let $V$ be $1$-dimensional simply connected oriented manifold and  let $A$ be an atlas of $\mathbb{R}$ on $V$ defining an orientation on $V$. Let $x_1$ and $x_2$ be two points of $V$ which are not separate. We will say that $x_1$ and $x_2$ are not separated on the right (resp.\ left) and write $x_1 \sim x_2$ mod $\lambda^{+}$ (resp.\ $\lambda^-$) if, given two charts $h_1$ and $h_2$ of $A$ such that $h_1(0)=x_1$ and $h_2(0)=x_2$, $h_2^{-1}h_1$ is defined for points with coordinates greater than $0$ (resp. less than $0$). These two relations are equivalence relations. Each class contains a finite or countable number of points. 

Let us call a foliation of $\mathbb{R}^2$ with consistent orientation on the leaves an oriented foliation and let us say two oriented foliations are equivalent if there is a self-homeomorphism of $\mathbb{R}^2$ that transports the oriented leaves of the first foliation to those of the second oriented foliation. (Notice that any foliation of $\mathbb{R}^2$ can be oriented.) We then show that the leaf space of $V$ of any oriented foliation of $\mathbb{R}^2$ is a $1$-dimensional oriented manifold whose order is canonically defined in each equivalence class modulo $\lambda^+$ or $\lambda^-$. We say that the oriented manifold $V$ is provided with an ordering. 

Two oriented foliations of $\mathbb{R}^2$ are equivalent if and only if there exists a homeomorphism of the leaf space of the first foliation to the leaf space of the second foliation that preserves orientation. 

Finally, to any simply connected $1$-dimensional manifold with countable basis that is orientable with an ordering, there corresponds a foliation of $\mathbb{R}^2$.

\begin{center}
    \includegraphics[width=130mm,scale=3.5]{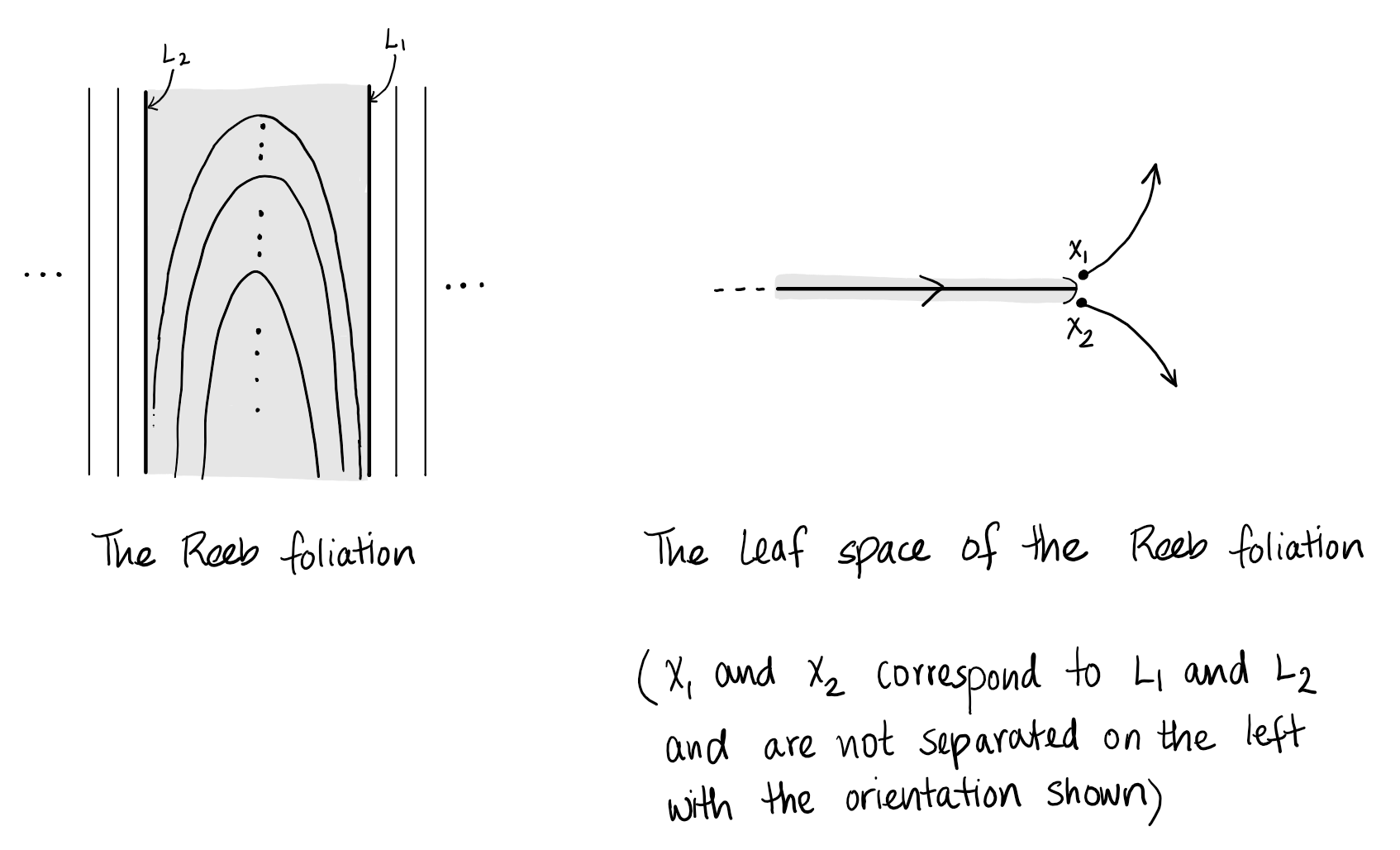}
\end{center}

\vspace{4mm}

\emph{Acknowledgments:} We wish to thank Andr\'e Haefliger for proofreading this translation. We also thank Indira Chatterji for helpful guidance and encouragement, and the referee for their detailed feedback and comments.

\end{document}